\documentclass[runningheads]{lncse}

\usepackage{amsmath}
\usepackage{amsfonts}
\usepackage{epsfig}
\usepackage{graphics} 

\DeclareFontFamily{OT1}{pzc}{}
\DeclareFontShape{OT1}{pzc}{m}{it}{<-> s * [1.2] pzcmi7t}{}
\DeclareMathAlphabet{\mathpzc}{OT1}{pzc}{m}{it}

\usepackage{bm}
\usepackage{amssymb}
\usepackage{mathtools}
\usepackage{stmaryrd}
\usepackage{tikz}
\usepackage{subcaption}
\usepackage{pgfplots}
\usepackage{pgfplotstable}

\numberwithin{equation}{section}

\newcommand{\Div}[0]{\mathbf{div}\,}
\newcommand{\bs}[1]{\boldsymbol{#1}}
\newcommand{\Sig}[0]{\bs{\sigma}}

\newcommand{\Eps}[0]{\bs{\varepsilon}}
\newcommand{\U}[0]{\bs{u}}
\newcommand{\V}[0]{\bs{v}}
\newcommand{\W}[0]{\bs{w}}

\newcommand{\F}[0]{\bs{f}}
\newcommand{\Z}[0]{\bs{0}}
\newcommand{\N}[0]{\bs{n}}

\newcommand{\T}[0]{\bs{t}}
\newcommand{\X}[0]{\bs{x}}

\newcommand{\jump}[1]{\left\llbracket #1 \right\rrbracket}

\newcommand{\dual}[1]{\left\langle #1 \right\rangle}
\newcommand{\Bf}[0]{\mathcal{B}}

\newcommand{\Sf}[0]{\mathcal{S}_h}
\newcommand{\vertiii}[1]{{\left\vert\kern-0.25ex\left\vert\kern-0.25ex\left\vert #1 
    \right\vert\kern-0.25ex\right\vert\kern-0.25ex\right\vert}}
    \newcommand{\enorm}[1]{\vertiii{#1}}

\newcommand{\Ch}[0]{\mathcal{C}_h}
\newcommand{\Eh}[0]{\mathcal{E}_h}
\newcommand{\Gh}[0]{\mathcal{G}_h}
\newcommand{\Nh}[0]{\mathcal{N}_h}

\newcommand{\Gcap}[0]{\mathcal{G}_h^{12}}

\newcommand{\un}{u_n}
\newcommand{\wn}{w_n}
\newcommand{\uhn}{u_{h,n}}
\newcommand{\vhn}{v_{h,n}}
\newcommand{\vn}{v_n}

\newcommand{\Q}{H^{\frac12}(\Gamma)}
\newcommand{\mQ}{H^{-\frac12}(\Gamma)}

 \newcommand{\situ}{\Sig_{i,t}(\U_i)} 

  \newcommand{\h}[0]{\mathpzc{h}}

 \newcommand{\bV}{\bm V}

\begin{document}

\title{Nitsche's Master-Slave  Method for \\
Elastic Contact Problems}

\titlerunning{Nitsche's master-slave method}

\author{ Tom Gustafsson\inst{1}  \and Rolf Stenberg\inst{2} \and Juha Videman\inst{3}   }

\authorrunning{Gustafsson et al.}   

\institute{
Department of Mathematics and Systems Analysis, Aalto University, 00076 Aalto, Finland
 {\tt tom.gustafsson@aalto.fi.}
  \and
 Department of Mathematics and Systems Analysis, Aalto University, 00076 Aalto, Finland
 {\tt rolf.stenberg@aalto.fi.}
 \and 
CAMGSD/Departamento de Matem\'atica, Universidade de Lisboa, Universidade de Lisboa, 1049-001 Lisbon, Portugal
  {\tt jvideman@math.tecnico.ulisboa.pt}
}

\maketitle

\begin{abstract} We survey the Nitsche's master-slave finite element method for elastic contact problems  analysed in 
\cite{GSV-elcontact}. The main steps of the error analysis are recalled and    numerical benchmark computations are presented.
\end{abstract}

\section{Introduction}

In a recent paper \cite{GSV-elcontact}, we studied Nitsche's method applied to contact problems between two elastic bodies. We considered three formulations, two of which take the different material properties of the  bodies into account. In the third method, which  will be detailed in this paper, the body with a~higher shear modulus is chosen as the master body and the slave one is mortared to it through Nitsche's method. We have the same error estimates for all three formulations but the master-slave approach appears to be the most straight-forward to implement.

Previously, the a priori estimates had been given under the assumption that the solution is in $H^s$, with $s>3/2$,  see {\sl e.g.} \cite{Frenchsurvey}, and the a posteriori estimates were derived using a saturation assumption. In \cite{GSV-elcontact}, we were able to improve the error analysis and avoid the saturation assumption. The key idea was to interpret Nitsche's method as a stabilised mixed method.

The plan of this paper is the following. In the next section we recall the elastic contact problem. In Section 3 we present the Nitsche's formulation, the stabilised method and show their equivalence. Then we summarise our error estimates and in the final section  give some numerical results supplementing those of \cite{GSV-elcontact}.

\section{The elastic contact problem}
\label{elastcont}
By $\Omega_i \subset \mathbb{R}^d$, $i=1,2$, $d=2,3$, we denote two elastic bodies in   contact, with the common boundary
 \mbox{$\Gamma = \partial \Omega_1 \cap \partial
\Omega_2$}.
 The parts of $\partial \Omega_i $  on which Dirichlet and Neumann boundary conditions are imposed are denoted by 
$ \Gamma_{D,i}$ and $ \Gamma_{N,i}$, respectively.
We let  $\U_i : \Omega_i \rightarrow
\mathbb{R}^d$   be the displacement of the body $\Omega_i$ and denote 
the strain tensors by 
$
    \Eps(\U_i) = \frac12 (\nabla \U_i + (\nabla \U_i)^T ).
$
The materials will be assumed to be isotropic and homogeneous, i.e., 
$ 
    \Sig_i(\U_i) = 2 \mu_i\,\Eps(\U_i) + \lambda_i\,\mathrm{tr}\,\Eps(\U_i) \boldsymbol{I},
$
where $\mu_i$ and $\lambda_i$ are  the  Lam\'e parameters.   We will exclude the possibility that the materials are nearly incompressible and hence it holds 
$\lambda_i \lesssim \mu_i$.  We assume thar $\mu_1\geq \mu_2$ and call the body $\Omega_1$ the \emph{master} and $\Omega_2$ the \emph{slave}. 
The outward unit normals to the boundaries are denoted by $\N_i  $ and we define $\N =  \N_1 =- \N_2$. Moreover, $\T$ denotes unit tangent vector satisfying $\N \cdot \T = 0$.

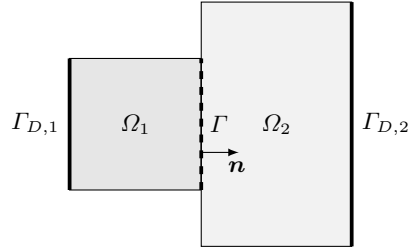
\begin{figure}[h!]
    \centering
    \begin{tikzpicture}[scale=0.5]
        \draw[fill=gray!20] (0.5, 0.5) rectangle (4, 4) node[pos=.5] {$\Omega_1$};
        \draw[fill=gray!10] (4, -1) rectangle (8, 5.5) node[pos=.5] {$\Omega_2$};
        \draw (4, 0.5) -- (4, 4) node[midway, anchor=west] {$\Gamma$};
        \draw[line width=0.5mm] (0.5, 0.5) -- (0.5, 4) node[midway, anchor=east] {$\Gamma_{D,1}$};
        \draw[line width=0.5mm] (8, -1) -- (8, 5.5) node[midway, anchor=west] {$\Gamma_{D,2}$};
        \draw[line width=0.5mm, dashed] (4, 0.5) -- (4, 4);
        \draw[-latex] (4, 1.5) -- (5, 1.5) node[pos=0.95, anchor=north] {$\boldsymbol{n}$};
    \end{tikzpicture}
    \caption{Notation for the elastic contact problem.}
    \label{fig:twobodynot}
\end{figure}
The traction vector    $\Sig_i(\U_i)\N_i $ is decomposed into its normal and tangential parts, viz.
\begin{equation}
  \Sig_i(\U_i)\N_i =\Sig_{i,n}(\U_i)+\Sig_{i,t}(\U_i).
  \end{equation}
For the scalar normal  tractions we  use the sign convention
\begin{equation}
\sigma_{1,n}(\U_1) = \Sig_{1,n}(\U_1)\cdot \N_1 ,
\  \mbox{ and } \  
\sigma_{2,n}(\U_2) = -\Sig_{2,n}(\U_2)\cdot \N_2,
\end{equation}
and note that on $\Gamma$ these tractions are either both zero or continuous and compressive, i.e. it holds that
\begin{equation}
\sigma_{1,n}(\U_1)=\sigma_{2,n}(\U_2), \quad \sigma_{i,n}(\U_i) \leq 0, \ i=1,2. 
\end{equation}

The physical non-penetration constraint on $\Gamma$  reads as 
$
      \U_1 \cdot \N_1 + \U_2 \cdot \N_2 \leq 0,
$
which, defining 
 $
\un= -( \U_1 \cdot \N_1 + \U_2 \cdot \N_2),
$
can be written as 
\begin{equation}
\jump{\un}\geq 0,
\end{equation}
where $\jump{\cdot }$ denotes the jump over $\Gamma$.

The \emph{Nitsche's method} is derived from the problem in displacement variables:
 find $\U=(\U_1,\U_2)$,   satisfying 
the \emph{equilibrium equations } for the two bodies
 \begin{equation}  \label{equlib}
  \ -\Div \Sig_i(\U_i)=\F_i \quad \text{in $\Omega_i$.}  
     \end{equation}
     The \emph{boundary conditions} are
   \begin{equation} 
        \U_i =\Z \quad   \mbox{on }  \Gamma_{D,i}, \quad 
        \Sig_i(\U_i) \N_i  =\Z \quad  \mbox{on }  \Gamma_{N,i} .
    \end{equation}    
    Next, we turn to  the common boundary. Here we assume that the \emph{tangential tractions  vanish}
    \begin{equation} 
        \situ = \Z \    \mbox{ on }\Gamma,
\end{equation}
    and that the \emph{normal stresses are continuous}
     \begin{equation} \label{conttrac}
        \situ = \Z, \quad
       \sigma_{1,n}(\U_1) - \sigma_{2,n}(\U_2)=0   \mbox{ on }\Gamma.
\end{equation}
The contact conditions are the \emph{non-penetration}
\begin{equation}
   \jump{\un}    \geq 0    \ \mbox{ on }\Gamma,
\end{equation}
and 
the \emph{non-positivity of the normal stresses}, and the \emph{compatibility condition}
\begin{equation}\label{non-p}
    \sigma_{i,n}(\U_i) \leq 0,  \quad 
          \jump{\un}  \sigma_{i,n}(\U_i)  = 0  \ \mbox{ on }\Gamma. 
\end{equation}
By the \emph{contact boundary} $\Gamma_C$ we mean the subset of $\Gamma$ wherein there is compression, i.e.~$ \sigma_{i,n}(\U_i) < 0  $.
On the complement $\Gamma \setminus \Gamma_C$, the normal tractions vanish. Note, however, that the contact  boundary is a priori unknown as it depends on the solution. 

The \emph{stabilised method} is based on the formulation in which the normal traction is an independent unknown. The equations \eqref{conttrac}  and 
 \eqref{non-p} are then replaced by
 \begin{align}
      \lambda + {\sigma_{1,n}(\U_1)}    =0,~~~\lambda +  {\sigma_{2,n}(\U_2)    }  &=0\  \text{ on } \Gamma, 
    \end{align}
and 
\begin{equation}
       \lambda \geq  0, \quad 
       \jump{\un} \lambda  = 0\   \mbox{ on } \Gamma. 
       \end{equation}

\section{The finite element methods}

The continuous displacements are in  $\bV = \bV_1 \times \bV_2$, with
$$
  \bV_i=\{\W_i   \in[H^1(\Omega_i)]^d : \W_i |_{\Gamma_{D,i}} = \Z\}. 
$$
By $\Ch^i$ we denote the  simplicial mesh on $\Omega_i$ which induces a facet mesh  $\Gh^i$ on $\Gamma$.  The finite element solution is sought in $ \bV_h = \bV_{1,h} \times \bV_{2,h}$, with 
$$
   \bV_{i,h} = \{ \V_{i,h} \in \bV_i : \V_{i,h}|_K \in [P_p(K)]^d~\forall K \in \Ch^i  \}.
$$

\smallskip
First, we recall the \emph{Nitsche's master-slave method}. 
To this end, we define the mesh function $\h_2$ on $\Gamma$ by 
$\h_2\vert_E=h_E$ for $E\in \Gh^2$. 
The method reads as follows:
  find $\U_h \in \bV_h$ such that
 \begin{eqnarray} \nonumber
         \sum_{i=1}^2 (\Sig_i(\U_{i,h}),\Eps(\V_{i,h}))_{\Omega_i}          
       &&  +\big\langle    \sigma_{2,n}(\U_{2,h})  ,  \jump{\vhn}   \big\rangle_{  \Gamma_c(\U_h)}+  
       \big\langle  \sigma_{2,n}(\V_{2,h})  ,  \jump{\uhn}   \big\rangle_{  \Gamma_c(\U_h)}  
            \\   
            && \hskip -50 pt
       +\gamma \big\langle   \frac{\mu_2}{\h_2} \jump{\uhn}  , \jump{\vhn}     \big\rangle_{  \Gamma_c(\U_h)}   
            = \sum_{i=1}^2 (\F_i,\V_{i,h})_{\Omega_i} \quad \forall \V_h \in \bV_h,
            \label{msn}
         \end{eqnarray}
 where $\gamma>C_I^{-1} $, with $C_I>0$ denoting the constant in the discrete trace inquality
    \begin{equation}
        C_I   h_E \| \sigma_{2,n}(\U_{2,h}) \|_{0,E}^2 \leq \mu_2\| \Sig_2(\U_{2,h})\|_{0,K}^2,  \quad E=K\cap \Gamma, 
            \label{dtrace}
    \end{equation}
    and 
 \begin{equation}
   \label{activecontact}
    \Gamma_c(\U_h) = \{ \X \in \Gamma : \sigma_{2,n}(\U_{2,h}) +\gamma \frac{\mu_2}{  \h_2} \jump{\uhn} < 0 \}.
\end{equation}
    The nonlinearity of the problem stems from this dependence of the contact boundary on the solution.

\smallskip

To define \emph{the stabilised method} we need some additional notation.
The normal traction $\lambda $ is in the space $H^{-\frac12}(\Gamma)$, dual to the trace space $H^{\frac12}(\Gamma)$, 
with  the norm $\Vert \cdot\Vert_{-\frac12, \Gamma}$ defined by duality. Defining
\begin{equation}
    \Bf(\W,\xi; \V,\eta)  = \sum_{i=1}^2 (\Sig_i(\W_i),\Eps(\V_i))_{\Omega_i} - \dual{\jump{\vn}  , \xi} - \dual{\jump{\wn} , \eta},
    \end{equation}
    and
    \begin{equation}\label{pospart}
    \varLambda=\{\xi \in \mQ : \dual{ w, \xi} \geq 0 ~~ \forall w \in \Q,~w \geq 0 ~\text{a.e.~on $\Gamma$}\},
\end{equation}
the mixed formulation of the problem is: find
 $(\U,\lambda) \in \bV \times \varLambda$ such that
    \begin{equation}
        \Bf(\U,\lambda; \V,\eta-\lambda) \leq \sum_{i=1}^2 (\F_i,\V_i)_{\Omega_i}  \quad \forall (\V,\eta) \in \bV \times \varLambda.
        \label{weakform}
    \end{equation}
The traction is approximated on the mesh $ \Gcap$ obtained as the intersection  
 of $\Gh^1$ and $\Gh^2$:
\begin{align}
  Q_h &= \{ \eta_h \in \mQ : \eta_h|_E \in P_p(E)~\forall E \in \Gcap  \}.
\end{align}
(Note that  since the approximation is discontinuous across element boundaries this is possible, even though the elements are general polygonals polyhedrons.)
Moreover, we introduce a subset of $\varLambda$,  denoted by $\varLambda_h $, as the positive part of $Q_h$, i.e.
\begin{equation}
    \varLambda_h = \{ \eta_h \in Q_h : \eta_h \geq 0 \} .
\end{equation}
The  stabilised bilinear form $\Bf_h$ is defined through
\begin{equation}
    \Bf_h(\W_h,\xi_h;\V_h,\eta_h) = \Bf(\W_h,\xi_h;\V_h,\eta_h) - \alpha \Sf(\W_h,\xi_h;\V_h,\eta_h),
\end{equation}
where $\alpha > 0$ is a stabilisation parameter and
\begin{equation}
    \Sf(\W_h,\xi_h;\V_h,\eta_h) =  \big\langle \frac{\h_2}{\mu_2} (\xi_h + \sigma_{2,n}(\W_{2,h})), \eta_h + \sigma_{2,n}(\V_{2,h})\big\rangle_\Gamma.
\end{equation}
The stabilised method is: find 
   $(\U_h,\lambda_h) \in \bV_h \times \varLambda_h$ such that
    \begin{equation}\label{VIh}
        \Bf_h(\U_h,\lambda_h; \V_h,\eta_h-\lambda_h) \leq \sum_{i=1}^2 (\F_i,\V_i)_{\Omega_i}\quad \forall (\V_h,\eta_h) \in \bV_h \times \varLambda_h.
    \end{equation}
 Note that
 \begin{equation}
    \Sf(\U_h,\lambda_h;\V_h,\eta_h) =  
    \big\langle \frac{\h_2}{\mu_2} (\lambda_h + \sigma_{2,n}(\U_{2,h})), \eta_h + \sigma_{2,n}(\V_{2,h})\big\rangle_\Gamma,
\end{equation}
 and hence the stabilised term amounts to a symmetric term including the residual $ \lambda_h + \sigma_{2,n}$.
 
 Now, by testing with $(\Z,\eta_h)$ in \eqref{VIh}, one can infer that
 
 \begin{equation}
\lambda_h=\big(
- \sigma_{2,n}(\U_{2,h}) - \frac{\mu_2}{ \alpha  \h_2} \jump{\uhn}\big)_+.
 \end{equation}
 Substituting this into the first equation obtained by testing with $(\V_h,0)$ in \eqref{VIh} we get the Nitsche's method \eqref{msn} with $\gamma=\alpha^{-1}$.

\section{Error estimates}

    The error estimate will be derived in the   norm
    \begin{equation}
    \enorm{(\W,\xi)}^2  = \sum_{i=1}^2\big( \mu_i \Vert \W \Vert_{1,\Omega_i} ^2+ \frac{1}{\mu_i} \| \xi \|_{-\frac12, \Gamma}^2\big).
    \end{equation} 
   The stability of the continuous problem  is given in the following theorem.
\begin{theorem}
    For every $(\W,\xi) \in \bV \times Q$ there exists $\V \in \bV$ such that
    \begin{equation}
        \Bf(\W,\xi;\V,-\xi) \gtrsim \enorm{(\W,\xi)}^2
 \ \mbox{ 
    and } \
        \|\V\|_V \lesssim \enorm{(\W,\xi)}.
    \end{equation}
\end{theorem}
 
The idea with stabilisation is that it yields a method which is always stable in a mesh-dependent norm for the Lagrange multiplier. Defining
\begin{equation}
\enorm{(\W_h,\xi_h)}_h^2=
\sum_{i=1}^2 \mu_i \Vert \W \Vert_{1,\Omega_i} ^2+\mu_2^{-1}\sum_{E\in \Gh^2}h_E \Vert \xi_h\Vert_{0,E}^2, 
    \end{equation}
we directly obtain the estimate.
\begin{theorem} Suppose that $0<\alpha < C_I$. Then,
    for every $(\W_h,\xi_h) \in \bV_h \times Q_h$, there exists $\V_h \in \bV_h$ such that
    \begin{equation}
        \Bf_h(\W_h,\xi_h;\V_h,-\xi_h) \gtrsim \enorm{(\W_h,\xi_h)}_h^2
 \ \mbox{   and }  \ 
        \|\V_h\|_V \lesssim \enorm{(\W_h,\xi_h)}_h.
    \end{equation}
\end{theorem}
 In view of Theorem 2,  the classical Verfürth trick yields the stability estimate in the correct norms.
 
 The error analysis then follows in a standard way, except for the additional term
 \begin{equation}
 \Big(\mu_2^{-1}\sum_{E\in \Gh^2}h_E \Vert \eta_h + \sigma_{2,n}(\V_{2,h}) \Vert_{0,E}^2\Big)^{1/2},
     \end{equation}
 where $(\V_h,\eta_h)$ are the interpolants of $(\U,\lambda).$
 However, by a posteriori error analysis techniques this term  can be bounded by
 \begin{equation}
 \enorm{(\U-\V_h,\lambda-\eta_h)} + \mathrm{HOT},
     \end{equation}
 where HOT stands for a higher order oscillation term.
 We thus arrive at the following quasi-optimality estimate of the method.

\begin{theorem} For $0<\alpha< C_I$ 
    it holds that
    \begin{eqnarray} 
        \enorm{(\U-\U_h, \lambda -   \lambda_h)}    \lesssim  &&  \inf_{(\V_h,\eta_h)\in  \bV_h\times\Lambda_h}\big(\enorm{(\U-\V_h,\lambda-\eta_h)}  +
           \sqrt{\langle \jump{u_n} ,\eta_h\rangle} \big)
           \nonumber 
            \\   
            && \qquad \qquad +
        \mathrm{HOT}.
         \end{eqnarray}
\end{theorem}

For the a posteriori error analysis,  we define the local estimators
\begin{alignat}{1}
    \label{apost1} \eta_K^2 &= \frac{h_K^2}{\mu_i} \|\hspace{0.2mm}\Div \Sig_i(\U_{i,h}) + \F_i \|_{0,K}^2, \quad K \in \Ch^i,
     \\
    \label{apost2} \eta_{E,\Omega}^2 &= \frac{h_E}{\mu_i} \left \| \jump{\Sig_i(\U_{i,h})\N} \right\|_{0,E}^2, \quad E \in \Eh^i, 
    \\
    \label{apost3} \eta_{E,\Gamma}^2 &= \frac{h_E}{\mu_i}  
    \|\hspace{0.2mm} \Sig_{i,t}(\U_{i,h})\|_{0,E}^2+
  \frac{\mu_i}{h_E} \Vert (\jump{u_{h,n} })_{-}\Vert_{0,E}^2,
    \qquad  E \in \Gh^i, 
                                     \\
    \label{apost4} \eta_{E,\Gamma_N}^2 &= \frac{h_E}{\mu_i} \left \|\hspace{0.2mm}\Sig_i(\U_{i,h})\N \right\|_{0,E}^2, \quad E \in \Nh^i,  
    \\
   \zeta_{E, \Gamma}^2&= \frac{h_E}{\mu_2} \left\| \lambda_h +\sigma_{2,n}(\U_{2,h}) \right\|_{0,E}^2,\quad E\in \Gh^2,
\end{alignat}
with $i=1,2$. The corresponding global estimator $\eta$ is then defined as
\begin{equation}
    \eta^2 = \sum_{i=1}^2 \Big\{\!\sum_{K \in \Ch^i} \eta_K^2 +\! \sum_{E \in \Eh^i} \eta_{E,\Omega}^2 +\! \sum_{E \in \Gh^i} \eta_{E,\Gamma}^2 +\! \sum_{E \in \Nh^i} \eta_{E,\Gamma_N}^2 \!\Big\}+\sum_{E\in \Gh^2}\zeta_{E,\Gamma}^2 \ .
\end{equation}

In addition, we need an estimator $S$ defined only globally as
\begin{equation}
    S^2 = \big\langle(\jump{u_{h,n}})_+,\lambda_h\big\rangle_\Gamma.
\end{equation}

\begin{theorem}[A posteriori error estimate]
    It holds that
    \begin{equation}
       \eta  \lesssim      \enorm{(\U-\U_h,\lambda-\lambda_h)} \lesssim \eta + S.
        \label{apostupper}
    \end{equation}
\end{theorem}

\section{Numerical experiments}

We investigate the performance of the master-slave method by solving adaptively
the problem \eqref{msn} using $P_2$ elements and the following geometry:
\begin{equation}
    \Omega_1 = [0.5,1] \times [0.25,0.75], \quad \Omega_2 = [1,1.6] \times [0,1].
\end{equation}
The boundary conditions are defined on
\begin{align}
    \Gamma_{D,1} &= \{ (x,y) \in \partial \Omega_1 : x=0.5 \}, \quad \Gamma_{N,1} = \partial\Omega_1 \setminus (\Gamma_{D,1} \cup \Gamma),\\
    \Gamma_{D,2} &= \{ (x,y) \in \partial \Omega_2 : x=1.6 \}, \quad \Gamma_{N,2} = \partial\Omega_2 \setminus (\Gamma_{D,2} \cup \Gamma),
\end{align}
while the material parameters are $E_1 = 1$, $E_2 = 0.1$ and $\nu_1 = \nu_2 = 0.3$.
The loading is
\begin{equation}
    \F_1 = (0,-\tfrac{1}{20}), \quad \F_2 = (0,0),
\end{equation}
which causes the left block to bend downwards. The active contact
boundary $\Gamma_c$ is sought by alternately evaluating
the inequality condition in \eqref{activecontact}
and solving the linearised problem with $\alpha = 10^{-2}$.

\begin{figure}
  \begin{subfigure}[t]{\textwidth}
    \centering
    \includegraphics[width=0.49\textwidth]{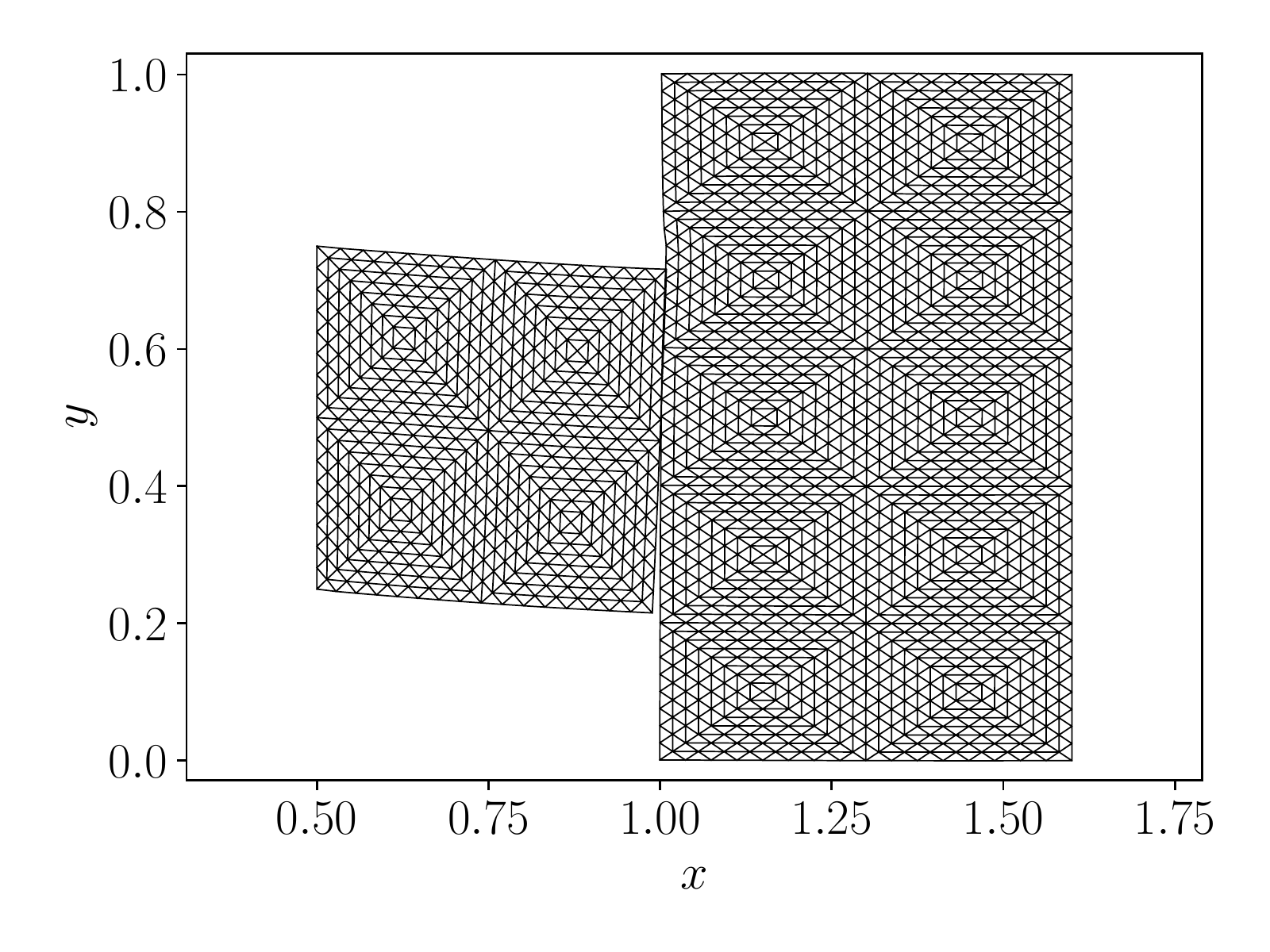}
    \includegraphics[width=0.49\textwidth]{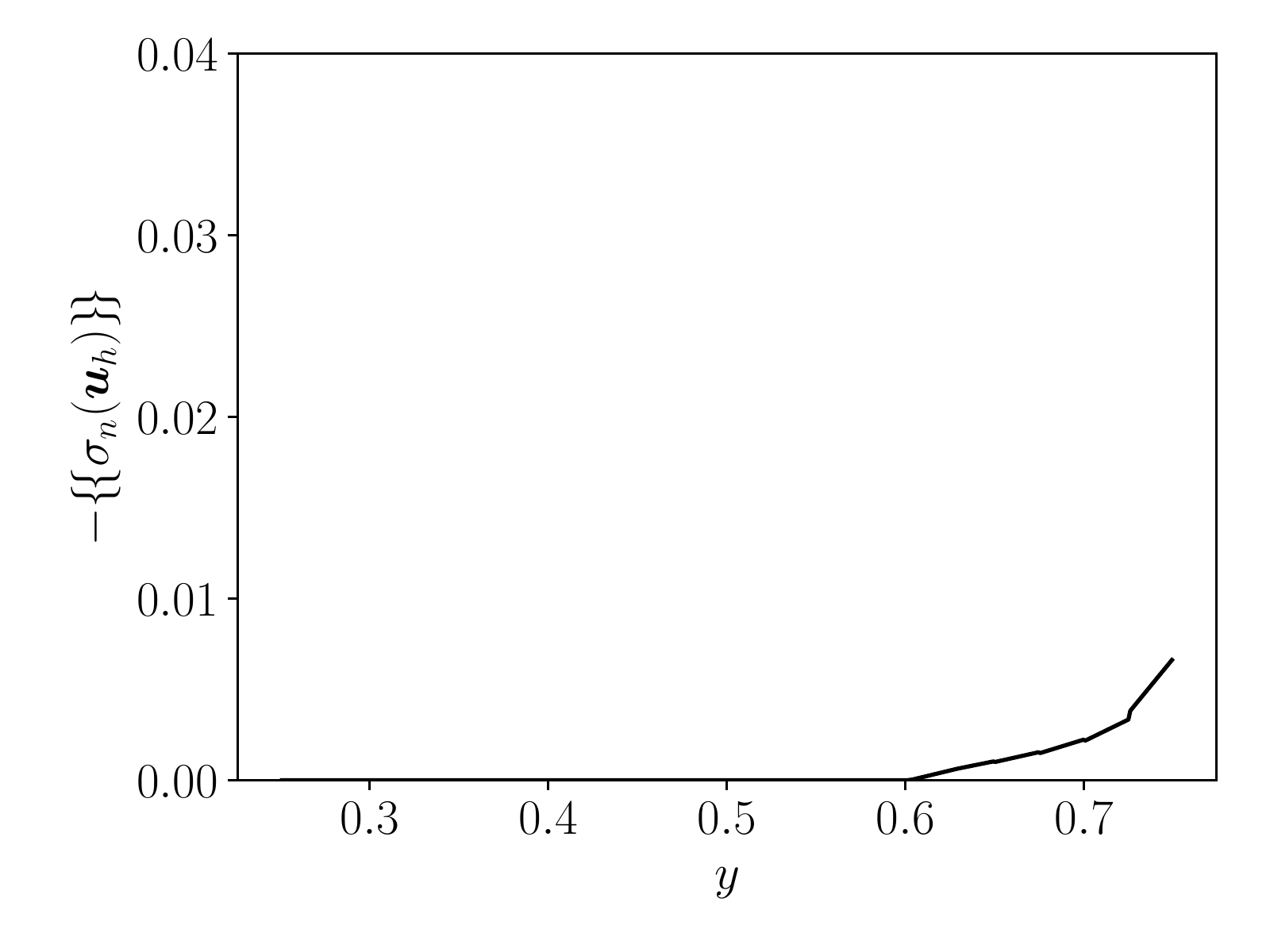}
    \caption{$P_2$ after 3 uniform refinements.}
    \label{fig:uniformseq}
  \end{subfigure}
  \\[0.2cm]
  \begin{subfigure}[t]{\textwidth}
    \includegraphics[width=0.49\textwidth]{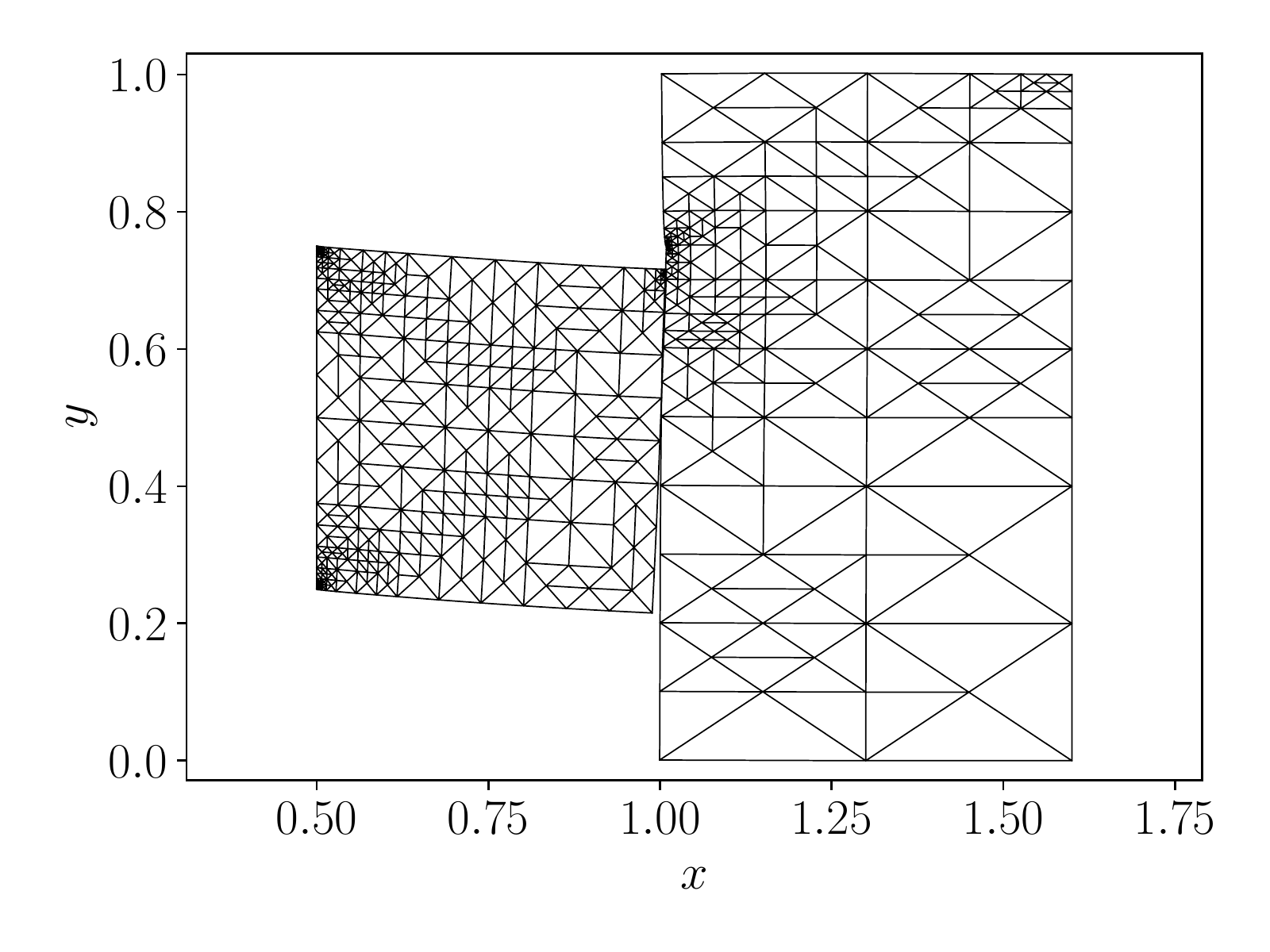}
    \includegraphics[width=0.49\textwidth]{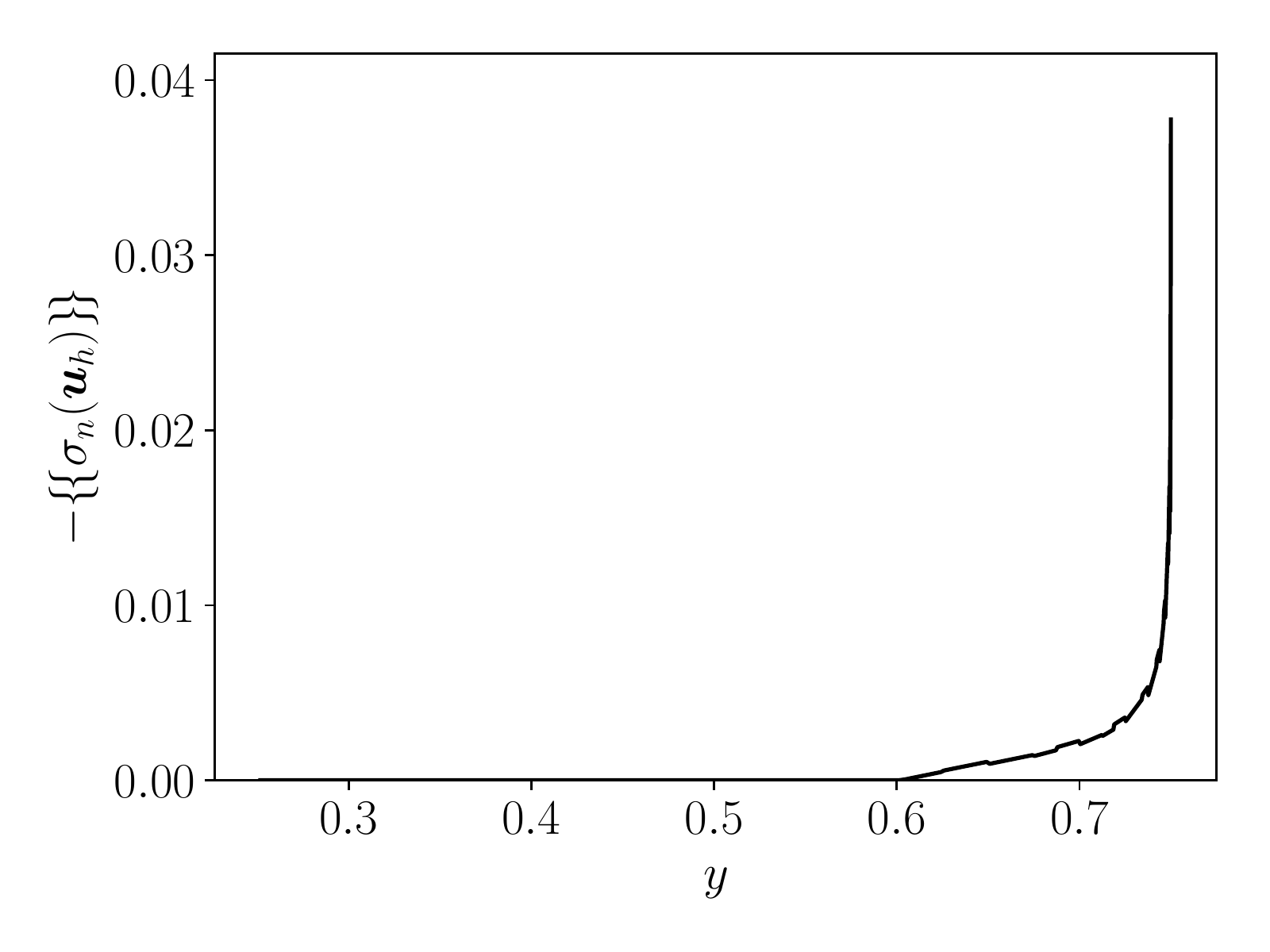}
    \caption{$P_2$ after 10 adaptive refinements.}
    \label{fig:adaptseq}
  \end{subfigure}
  \\[0.5cm]
  \begin{subfigure}[t]{\textwidth}
    \centering
    \begin{tikzpicture}[scale=1.0]
      \begin{axis}[
          xmode = log,
          ymode = log,
          xlabel = {$N$},
          ylabel = {$\eta+S$},
          grid = both,
          legend style={at={( 0.97,0.97)}, anchor=north east},
        ]
        \addplot table[only marks,x=ndofs,y=eta] {\expfirst};
        \addplot table[only marks,x=ndofs,y=eta] {\expsecond};

        \addplot[blue] table[y={create col/linear regression={y=eta}}] {\expfirst};
        \xdef\expfirstcoeff{\pgfplotstableregressiona};
        
        \addplot[red] table[y={create col/linear regression={y=eta}}] {\expsecond};
        \xdef\expsecondcoeff{\pgfplotstableregressiona};

        \addlegendentry{Uniform $P_2$, $\mathcal{O}(N^{\pgfmathprintnumber{\expfirstcoeff}})$}
        \addlegendentry{Adaptive $P_2$, $\mathcal{O}(N^{\pgfmathprintnumber{\expsecondcoeff}})$}
      \end{axis}
    \end{tikzpicture}
    \caption{The global a posteriori error estimator $\eta + S$ as a function of the number of degrees-of-freedom $N$.}
    \label{fig:convergence}
  \end{subfigure}
  \caption{The numerical example, the resulting meshes and global error estimators.}
  \label{fig:example}
\end{figure}
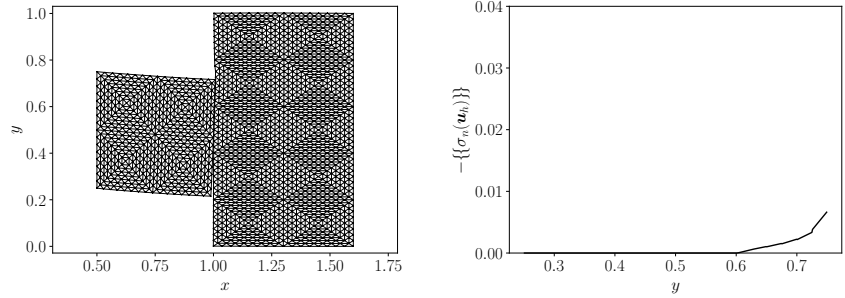
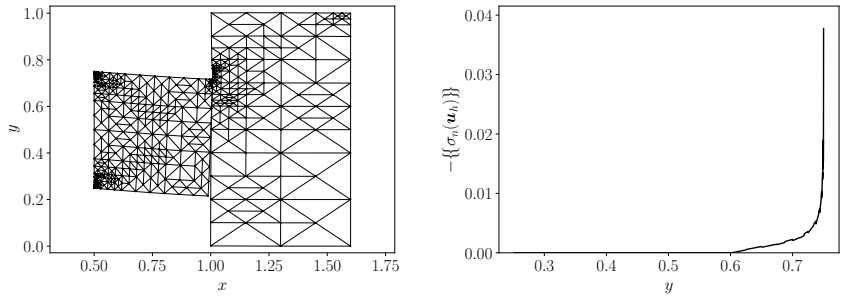
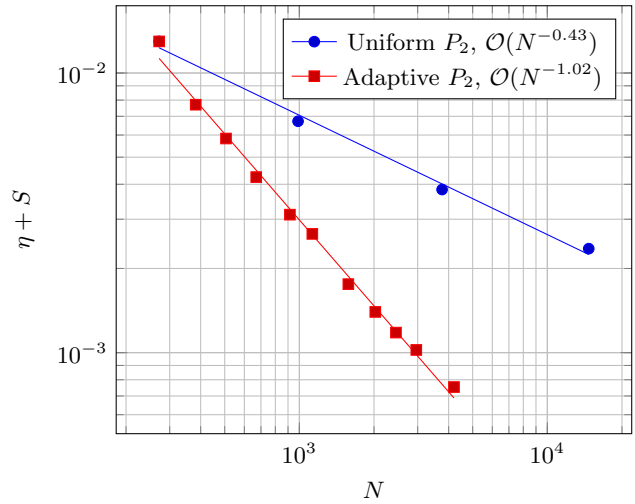

The final meshes and the respective approximation of the contact
force is given in Figure~\ref{fig:uniformseq} and \ref{fig:adaptseq} where we use the notation
$\{\!\{\sigma_n(\bs{u}_h)\}\!\}$ for the mean normal stress over the contact boundary.
The~resulting global a posteriori error estimator is given as a function of the number of degrees-of-freedom 
in Figure~\ref{fig:convergence}.
We observe, in particular, that the asymptotic rate of convergence for the
total error estimator is improved from $\mathcal{O}(N^{-0.43})$ to $\mathcal{O}(N^{-1.02})$
where the latter corresponds to the rate of convergence one expects from $P_2$ elements
and a completely smooth solution.

\pgfplotstableread{
  ndofs eta
  272 0.012971602162436794
  382 0.007703923242312411
  506 0.005831872626842647
  670 0.004244223838830169
  916 0.0031152584647078335
  1128 0.0026587938510649425
  1576 0.001758130298760318
  2024 0.0013984885158793217
  2452 0.0011796605044011236
  2962 0.0010226639255012649
  4204 0.0007534454969743955
}\expsecond 

\pgfplotstableread{
  ndofs eta
  272 0.012971602162436794
  988 0.006722333249642182
  3764 0.0038307602868424863
  14692 0.0023545343491632965
}\expfirst

\section*{Acknowledgments}  
This work was (partially) funded by the Portuguese government through FCT - Funda\c{c}\~ao para a Ci\^encia e a Tecnologia, I.P., under the projects PTDC/MAT-PUR/28686/2017, and UID/MAT/04459/2013,  and by the \break Academy of Finland, decision nr. 324611.
%
%
%

\ifx\undefined\bysame
\newcommand{\bysame}{\leavevmode\hbox to3em{\hrulefill}\,}
\fi

\bibliographystyle{siam}
\bibliography{nitsche_mortaring,ps-ref}

\end{document}